\documentclass{amsart}
\usepackage{amsmath,amsfonts,amssymb,hyperref}

\newcommand{\al}{\alpha}    \newcommand{\be}{\beta}
\newcommand{\de}{\delta}    
\newcommand{\ep}{\epsilon}  
\newcommand{\la}{\lambda}   
    
\newcommand{\om}{\omega}    \newcommand{\Om}{\Omega}
\newcommand{\ga}{\gamma}    
\def\<{\langle}             \def\>{\rangle}
\newcommand{\R}{\mathbb{R}}

\newcommand{\Sp}{\mathbb{S}}

\newcommand{\Sc}{\mathcal{S}}
\newcommand{\g}{{\mathrm g}}
\newcommand{\hf}{\frac 12}

\newcommand{\pt}{\partial_t}\newcommand{\pa}{\partial}
\newcommand{\les}{{\lesssim}}

\newcommand{\beeq}{\begin{equation}}\newcommand{\eneq}{\end{equation}}

\theoremstyle{plain}
\newtheorem{thm}{Theorem}

\newtheorem{lem}[thm]{Lemma}

\theoremstyle{remark}
\newtheorem{rem}{Remark}

\theoremstyle{definition}

\newenvironment{prf}{\noindent {\bf Proof.} }{\endprf\par}
\def \endprf{\hfill  {\vrule height6pt width6pt depth0pt}\medskip}

\numberwithin{equation}{section}

\begin{document}

\title[Strauss Conjecture]{Recent works on the Strauss conjecture}
\thanks{The first author was supported by the Fundamental Research Funds for the Central Universities, NSFC 10871175 and
10911120383.}

\author{Chengbo Wang}
\address{Department of Mathematics, Zhejiang University, Hangzhou 310027,
China} \email{wangcbo@gmail.com}\urladdr{http://www.math.zju.edu.cn/wang}

\author{Xin Yu}
\address{EPGY, Stanford University, Stanford, CA 94305,
USA} \email{yumei165@stanford.edu}

\date{}
\dedicatory{} \commby{}

\maketitle \tableofcontents

\section{Introduction}
There have been many exciting breakthroughs on the understanding of the long time existence v.s. blow up for the semilinear wave equations with small initial data, which is also known as the Strauss conjecture. In this paper, we would like to present some highlights occurred in these works.

Most of these breakthroughs occurred in the context of the problem with nontrapping obstacles. However, to prevent technical complexities in dealing with obstacles, we will mainly discuss the problem without obstacles and illustrate the methods to deal with the problem with obstacles.

Let $\Box=\pt^2-\Delta$,
we shall consider wave equations in the Minkowski space $\R_+\times \R^n$
\begin{equation}\label{LW}
\begin{cases}
\Box u(t,x)=F_p(u(t,x)), \quad (t,x)\in \R_+\times \R^n
\\
u(0,\cdot)=f, \quad
\partial_t u(0,\cdot)=g\ .
\end{cases}
\end{equation}
We shall assume that the nonlinear term behaves like $|u|^p$ when $u$ is small, and so we assume that
\begin{equation}\label{Fp-assume}
\sum_{0\le j\le k} |u|^j\, \bigl|\, \partial^j_u F_p(u)\, \bigr| \,
\lesssim \, |u|^p
\end{equation}
when $u$ is small, for some $k$ to be fixed later. We shall also assume throughout
that the spatial dimension satisfies $n\ge 2$. We will denote $|D|=\sqrt{-\Delta}$ as the square root of the positive Laplacian.

\subsection{Initial analysis of the equation: scaling}
Let us consider the model equation
\beeq\label{slw-model}\Box u=|u|^p\ .\eneq
It is easy to check that the equation is invariant under the scale transfrom
$$u(t,x)\rightarrow u_\la(t,x)\equiv \la^{-2/(p-1)}u(t/\la,x/\la),\ \la>0,$$
that is, if $u(t,x)$ is a solution to the equation, $u_\la$ is also a solution for any $\la>0$.

It is known that the scaling invariance will give us a lower bound of the Sobolev regularity of the well posedness for the problem (see e.g. Section 3.1 of \cite{Tao06} for general discussion and \cite{FaWa05}, \cite{FaWa08} for the results of ill-posedness for the problem with low regularity).

Recall that the homogeneous Sobolev space $\dot H^s$ with $s<n/2$ is defined as the completion of $C_0^\infty$ in the space of temper distribution $\Sc'$ under the semi-norm
$$\|f\|_{\dot H^s}=\|(-\Delta)^{s/2}f\|_{L^2}\ .$$
The scaling invariance gives us the scaling regularity \beeq\label{sc}s_c=\frac n2-\frac 2{p-1}\ .\eneq

\subsection{Strauss conjecture: the beginning of the story}
Let us now give a bit of historical background. In 1979, Fritz John \cite{John79} showed that
when $n = 3$ global solutions always exist for \eqref{slw-model} if $p > 1 + \sqrt{2}$ and the size of the compactly supported smooth initial data, $\ep > 0$, is small. He
also showed that the power $1 +\sqrt{2}$ is critical in the sense that no such result can
hold if $p < 1 +\sqrt{2}$.

The number $1 + \sqrt{2}$ appears first in Strauss' work \cite{Strauss78} on
scattering for small-amplitude semilinear Schr\"{o}dinger equations.
Based on this, he made the insightful conjecture in \cite{Strauss81} that
the critical power in John's theorem in $n\ge 2$ space dimensions, denoted by $p_c(n)$, is
the positive root of the quadratic equation $$(n - 1)p^2 - (n + 1)p - 2 = 0\ .$$ In other words, \beeq\label{pc}p_c(n)=\frac{n+1+\sqrt{(n+1)^2+8(n-1)}}{2(n-1)}=\frac{n+1+\sqrt{n^2+10n-7}}{2(n-1)}\ .\eneq
In particular,
$$p_c(2)=\frac{3+\sqrt{17}}2, \ p_c(3)=1+\sqrt{2},\ p_c(4)=2\ .$$
As a side remark, for the case $n=1$, it is easy to see that we can not have the global existence with small data for any $1<p<\infty$, see e.g. Kato \cite{Kato80}.

Basically, there are two parts in the Strauss conjecture: global existence with small data for $p>p_c$ and finite in time blow up for $p<p_c$. The main line of research is about the part of global existence with small data for $p>p_c$.

Before the discussion of the history of resolving the Strauss conjecture, let us record an interesting observation about the relations between the scaling regularity $s_c$ and the critical power $p_c$. Besides the critical regularity $s_c$, there is one more Sobolev regularity, namely \beeq\label{sd}s_d=1/2-1/p,\eneq as far as the radially symmetric functions are concerned (see Sogge \cite{So08} Section IV.4). It can be shown that $s_d$ is a lower bound of the Sobolev regularity for local existence results when $n=2, 3$, even for the radial data (Hidano \cite{Hidano-p}). It is interesting to note that, for $p>0$,
\beeq\label{observ1}s_c>s_d \Leftrightarrow p>p_c\ .\eneq
As a side remark, there is similar phenomenon for the equation $\Box u=|\pt u|^p$ with $3/2$ as the role of $s_d$, which is also known as the Glassey conjecture (see \cite{HWY2}).

We would like to categorize the history of the study of the Strauss conjecture into roughly three periods.
\subsection{First wave of study: 1979-1985}
Shortly after the pioneering work of John \cite{John79}, Glassey extends the result to the case $n=2$ in \cite{G}, \cite{Glassey}. John's
blow up result was then extended by Sideris \cite{Sideris}, showing that, for all $n\ge 4$, the solution can  blow up for arbitrarily small data if $1<p < p_c$. Sideris' proof of the blow up result is quite delicate, using sophisticated computation involving spherical harmonics and other special functions. His proof was simplified by Jiao and Zhou \cite{JiaoZhou03} and Yordanov and Zhang \cite{YoZh}.

The critical case $p=p_c$ was not clear until the work of Schaeffer \cite{Schaeffer}, where blowup was shown when $n=2, 3$.

Here is a comparison of the results in the first wave of study:

\begin{tabular}{|l|l|l|l|}
\hline & $n=2$ & $n=3$ & $n\ge 4$\\\hline $1<p<p_c$ blow up &
Glassey 81 & John 79 & Sideris 84\\\hline $p>p_c$ global  & Glassey
81    & John 79 & open \\\hline
 $p=p_c$ blow up & Schaeffer 85    & Schaeffer 85 & open  \\\hline
\end{tabular}

For the part of global existence, the authors mainly exploit the positivity of the fundamental solution of the wave equation when $n=2,3$ (together with the Huygens
principle). More precisely, they prove the global existence by iteration in the space of the form
$$\|\left<t\right>^{\frac{n-1}2}\left<t-|x|\right>^{q(p)} u(t,x)\|_{L^\infty_t L^\infty_x}$$
for some $q(p)$.

\subsection{Second wave of study: 1990-1997}
The main themes in the second wave of study (mainly 1990-1997) are to obtain more precise information about the lifespan (time of existence), i.e., to obtain the sharp life span (denoted by $T_\ep$) for $p\le p_c$ and prove the global existence for high dimensions.

Before discussing the history in this period, let us discuss the expected sharp lifespan, denoted by $L_\ep$, for the equation with data of size $\ep>0$. If $p>p_c$, we certainly have $L_\ep(p)= \infty$. If $1<p<p_c$, the sharp lifespan is expected to be
of order
\beeq\label{lifespan-le}L_\ep(p)\simeq \ep^{2 p (p-1)/((n-1)p^2-(n+1)p-2)}\ .\eneq
Here by $A\simeq B$ we mean that there exist absolute constants $C, c>0$ such that $c B < A < C B$. Else, if $p=p_c$, we expect to have
\beeq\label{lifespan-cri}\log L_\ep(p)\simeq
\ep^{-p (p-1)}\ .\eneq

Here, we record one more interesting observation about the relation between the scaling regularities $s_c$, $s_d$ and the lifespan for $p<p_c$:
\beeq\label{observ2}L_\ep(p)\simeq \ep^{1/(s_c-s_d)}\ .\eneq

For the part of global existence for $n\ge 4$ and $p>p_c$, Zhou proved the case $n=4$ in \cite{Zh95}. Kubo \cite{Kubo96}
has verified it under the assumption of spherical symmetry in odd spatial dimensions.
He also announced a proof of the corresponding result in even dimensions.
Shortly after, Lindblad and Sogge \cite{LdSo96} proved the global existence for all dimensions if one assumes radial symmetry and when $n \le 8$ if this assumption is removed. Finally, Georgiev, Lindblad and Sogge \cite{GLS97} proved the result for all dimensions $n\ge 4$. Later, a simple proof was given by Tataru \cite{Ta01-2}.

On the other hand, for the estimate of the sharp lifespan for $p\le p_c$, the first results were obtained for $n\le 3$ in the works of Lindblad \cite{Ld90} and Zhou \cite{Zh92}, \cite{Zh92-2} \cite{Zh93} (see also  Takamura \cite{Takamura}). When $n=4$ and $p=p_c$, Li and Zhou \cite{LiZhou95} proved almost global existence with $T_\ep\ge L_\ep$. Additionally, Lindblad and Sogge
\cite{LdSo96} obtained optimal lifespan estimates $T_\ep\ge L_\ep$ for all subcritical
powers $1 < p \le p_c$ if $n \le 8$, while in higher dimensions they obtained such results
when $1 < p \le (n + 1)/(n - 1)$ or the initial data are radial symmetric. The upper bound of the lifespan $T_\ep\le L_\ep$ was verified in Zhou-Han \cite{ZhouHan} for $n\ge 3$ and $p<p_c$.

The problem of blow up for $p=p_c$ and $n\ge 4$ was unsolved until the works of Yordanov-Zhang \cite{YorZh06} and Zhou \cite{Zh07}. The upper bound estimates of the lifespan $T_\ep\le L_\ep$ for $n\ge 3$ and $p=p_c$ were obtained very recently in the works of Takamura-Wakasa \cite{TaWa} and Zhou-Han \cite{ZhouHan2}. The expected lower bound estimates of the lifespan $T_\ep\ge L_\ep$ are still open for $n\ge 9$ and $p\le p_c$, which is the last open problem in the Strauss conjecture (a weaker lower bound for $p=p_c$ was obtained by Tataru \cite{Ta01-2}).

Here is a comparison of the results in the second wave of study:

\begin{tabular}{|l|l|l|l|}
\hline & $n=2, 3$ & $4\le n\le 8$ & $n\ge 9$\\\hline $p<p_c$ & sharp
(\cite{Ld90}, \cite{Zh92}, \cite{Zh93}) & sharp (\cite{LdSo96}, \cite{ZhouHan}) & radial sharp (\cite{LdSo96}, \cite{ZhouHan})\\\hline $p>p_c$ & 1980's
& global (\cite{Zh95}, \cite{LdSo96}) & global (\cite{GLS97}, \cite{Ta01-2})
\\\hline
 $p=p_c$ & sharp (\cite{Zh92}, \cite{Zh93}) & sharp & lower bound (\cite{Ta01-2}) \\
  &&(\cite{LdSo96}, \cite{TaWa}, \cite{ZhouHan2}) & upper bound (\cite{TaWa}, \cite{ZhouHan2}) \\\hline
\end{tabular}

For the proof of global existence, Lindblad and Sogge \cite{LdSo96} proved the result by iteration in the space
$$\||x|^{(3-n)/2}u(t,x)\|_{L^{pq}_t L^p_{|x|} L^{2n/(n-1)}_\theta}\ ,$$
for $q=q(p)$ so that (with the same scaling as $\|u\|_{\dot H^{s_c}}$)
$$\frac 1q=\frac 2{p-1}-\frac{n-1}2\ .$$
Here we are using the mixed-norm notation
\beeq\label{label}\|u\|_{L^q_t L^r_{|x|} L^\sigma_\theta}^q=
\int_\R\left(\int_0^\infty\left(\int_{\Sp^{n-1}} \left|u(t,
|x|\omega)\right|^\sigma d\omega\right)^{r/\sigma}|x|^{n-1} d|x|
\right)^{q/r} dt\eneq (so that it is the same as $L^q_t L^r_x$ if $r=\sigma$).

Instead, in \cite{GLS97}, the authors exploited the estimates for the inhomogeneous wave equations in the weighted Lebesgue space
$$\|(1+|t^2-|x|^2|)^{q} u\|_{L^{p+1}_{t,x}(\R_+\times \R^{n+1})}$$
for some $q(p)$ such that
$$\frac{1}{p(p+1)}<q(p)<\frac{n-1}2-\frac n{p+1}\ .$$

\subsection{Third wave of study: 2008-now}
The classical Strauss conjecture (i.e., the problem on the flat Minkowski spacetime) has been solved, except the problem of existence for $n\ge 9$ and $p\le p_c$ with the expected lower bound for the lifespan $T_\ep\ge L_\ep$. It is still interesting, however, to investigate the problem on various manifolds. Recently, there have been many exciting breakthroughs on the understanding of this problem on the various manifolds, in particular the space with nontrapping obstacles.

We will concentrate on the case of the problem with nontrapping obstacles, since most of the new ideas were developed in this setting. The other examples of the manifolds include the asymptotically flat/Euclidean manifolds and the Schwarzschild/Kerr spacetime. For the results on the asymptotically Euclidean manifolds, we refer the interested reader to the recent works of Sogge and Wang \cite{SoWa10} and our paper \cite{WaYu11}, where the global existence for $p>p_c$ and $n=3,4$ was proven in the case of the short range perturbation of the flat metric.
For the results on the Schwarzschild spacetime,  we refer the interested reader to the works of Dafermos and Rodnianski \cite{DaRo05} and Blue-Sterbenz \cite{BlSt06}, where the global existence for $p>3$ was proven. Recently, Lindblad, Sogge and the first author are able to prove the result for $p>2.5$, by essentially adapting the method of \cite{DMSZ} and exploit the local energy estimates of \cite{MMTT10}. We expect that the sharp result $p>p_c$ can be proven in the near future.

To handle the problem with nontrapping obstacles, we typically need to develop new ideas/methods of proof for the problem without the obstacles and then try to adapt the proof in the new setting.

The first results in this direction were obtained in the works of Du and Zhou \cite{DuZh08} and Du, Metcalfe, Sogge and Zhou \cite{DMSZ}. The common feature in these two papers is to exploit the weighted space-time $L^2$ estimates (which is known as local energy estimates, KSS estimate \cite{KeSmSo02} or Morawetz-KSS estimates) with respect to the norm
$$\|\<x\>^{-\al} u\|_{L^2_{t,x}}$$
for certain choices of $\al$. In \cite{DuZh08}, the authors studied the Dirichlet problem of quasi-linear wave equation (involving the nonlinearity $u^2$) outside of star-shaped obstacle in three space dimensions, obtaining the sharp lower bounds of the lifespan $T_\ep\ge L_\ep$. In \cite{DMSZ}, the authors verified the part of global existence for $p>p_c$, $n=4$ and the Dirichlet problem with nontrapping compact obstacles. In addition, they can prove almost global existence with a weak lower bound (comparing to $L_\ep$) $$T_\ep\ge A_\ep:=\exp(c/\ep)$$ for $p=p_c=2$.

On the other hand, Hidano \cite{Hi07} proposed a new proof of the Strauss conjecture when $n=2,3,4$ and the initial data are radial, where he proved weighted Strichartz estimates of the type
$$\||x|^{-\al} u\|_{L^p_{t,x}}$$ for the wave equations with radial data.
Hidano's idea was further developed in the works of Fang and Wang \cite{FaWa11} and Hidano, Metcalfe, Smith, Sogge and Zhou \cite{HMSSZ}, where the weighted Strichartz estimates of the type $$\||x|^{-\al} u\|_{L^p_{t,|x|}L^2_\theta}$$ were proven and so is the Strauss conjecture for $n=2,3,4$. In \cite{FaWa11}, the weighted Strichartz estimates with angular smoothing effect were proven for the wave equation and the Schr\"{o}dinger equation as well. In \cite{HMSSZ}, the results were also generalized to the Dirichlet and Neumann problem with nontrapping obstacle for $n=3,4$.
As was shown in \cite{HMSSZ}, though, a limitation arises in the
proof which is only relevant when the spatial dimension, $n$, equals two.
This is because that the $TT^*$ arguments involving the Christ-Kiselev lemma \cite{ChKi01} a priori require that the Sobolev regularity for the data in the homogeneous estimates be equal to $1/2$ when $n = 2$, with similar restrictions on the estimates for the inhomogeneous wave equation. The result of \cite{HMSSZ} was generalized to certain trapping obstacles in the second author \cite{Yu11}, under the condition that there is local energy estimates with loss of regularity.

The difficulty for $n=2$ was remedied recently in the work of Smith, Sogge and the first author \cite{SmSoWa10p}, where, instead of proving the weighted Strichartz estimates, the authors proved a generalized version of the Strichartz estimates without weights (we will call it generalized Strichartz estimates)
$$\|u\|_{L^q_t L^r_{|x|} L^2_\theta}\ .$$

On the other hand,  the upper bound of the lifespan $T_\ep\le L_\ep$ was verified in Zhou-Han \cite{ZhouHan} for $n\ge 3$ and $p<p_c$, for the problem with Dirichlet boundary condition.

Here is a comparison of the results in the third wave of study:

\begin{tabular}{|l|l|l|l|l|}
\hline & $n=2$& $n=3$ & $n=4$& $n\ge 5$\\\hline $p<p_c$ & open & \cite{Yu11}  $T_\ep\ge L_\ep$, \cite{ZhouHan} $T_\ep\le L_\ep$ & \cite{ZhouHan} $T_\ep\le L_\ep$ & \cite{ZhouHan} $T_\ep\le L_\ep$
\\\hline $p>p_c$ & \cite{SmSoWa10p}
& \cite{HMSSZ} & \cite{DMSZ}, \cite{HMSSZ} & open
\\\hline
 $p=p_c$ & open & \cite{Yu11}  $T_\ep\ge A_\ep$ &
 \cite{DMSZ}  $T_\ep\ge A_\ep$ & open
  \\\hline
\end{tabular}

\subsection{Structure of the paper} We will present and prove these three essential estimates (Morawetz-KSS estimates, weighted Strichartz estimates and generalized Strichartz estimates), give the proof of the Strauss conjecture without obstacles, and illustrate the method to deal with the problem with obstacles. In the appendix, we give the proofs for some basic estimates, including the trace lemma, Morawetz type estimates (also known as local energy estimates), and Morawetz-KSS estimates,  which are also of independent interest.

\section[Morawetz-KSS estimates]{Morawetz-KSS estimates and the global existence for $n=4$}

\subsection{A Proof of the Strauss conjecture: $n=4$}

Intuitively, consider the equation $\Box u=|u|^p$, we can use Morawetz type estimate \eqref{8} to provide a proof of the global existence for $n=4$. The idea is that if we write $|u|^p=|u|^2 |u|^{p-2}$, the term $|u|^{p-2}$ will provide certain spatial decay at infinity (since $p>2$), which can be absorbed in the term $|u|^2$ to go back to the local energy norm. More precisely, recall the duality of the trace estimate \eqref{7-g-Trace} with $s=1$ gives us
$$\|v\|_{\dot H^{-1}}\les \||x|^{-(n-2)/2}v\|_{L^1_{|x|}L^2_\theta}\ ,$$
we can estimate the solution $u$, with initial data of size $\ep$, as follows,
\begin{eqnarray*}
  \|\<x\>^{-1/2-\de}u\|_{L^2_{t,x}}&\les& \|f\|_{L^2_x}+\|g\|_{\dot H^{-1}}+\||u|^p\|_{L^1_t \dot H^{-1}}\\
  &\les& \ep+\||x|^{-(n-2)/2}|u|^p\|_{L^1_t L^1_{|x|} L^2_\theta}\\
  &\les& \ep+\||x|^{-(n-2)/2-2\de}|u|^2\|_{L^1_t L^1_{|x|} L^2_\theta}\les \ep+LHS^2\ ,
\end{eqnarray*} which will provide a global bound of the local energy norm. Here, only the last two steps are not rigorous.

To give a rigorous proof of the global existence for $n=4$, we need to use the spatial decay estimate for $R>1$,
\beeq\label{2-decay}
\Vert h \Vert_{L^{\infty} (R \leq \vert x \vert \leq R+1)}
\lesssim R^{- 3/2} \sum_{\vert \alpha \vert \leq 2, |\be|\le 1}
\Vert \Omega^{\alpha}\pa_x^\beta h \Vert_{L^2(R-1 \leq \vert x \vert
\leq R+2)}\ ,\eneq
where $\Omega$ denotes the rotation vector field. This inequality follows from the trace estimate \eqref{eq-trace} and the Sobolev lemma on the sphere $H^2_\omega(\Sp^3)\subset L^\infty_\omega$. Also the trace estimate \eqref{7-g-Trace} with $s=1$ together with the Sobolev embedding $\dot H^1\subset L^{4}(\R^4)$ yield
\begin{eqnarray}
  \|v\|_{\dot H^{-1}}&\les& \||x|^{-1}v\|_{L^1_{|x|}L^2_\omega(|x|\ge 1)}+
\|v\|_{L^{4/3}(|x|<1)}\nonumber\\
&\les&
\|\<x\>^{-1}v\|_{L^1_{|x|}L^2_\omega}+
\|v\|_{L^2}\ .\label{2-decay2}\end{eqnarray}

Now let us solve the equation by iteration. Set $u_0\equiv 0$, and define $u_k$, $k=1,2,3,\dots$ inductively as follows
$$\Box u_k=F_p(u_{k-1}),$$
with initial data $u(0,\cdot)=f$ and $\partial_t u(0,\cdot)=g$. We need to show the sequence is uniformly bounded and Cauchy in some sense. Here we just give the proof of the uniform boundedness of
$$M_k=\sum_{|\al|\le 2, |\be|\le 1}\left(\|\<x\>^{-1/2-\de}Z^\al\pa^\be u_k\|_{L^2_{t,x}}+\|Z^\al \pa^\be u_k\|_{L^\infty_t L^2_{x}}\right),$$
where $Z=\{\pa, \Omega\}$ with $\pa=(\pa_t, \pa_x)$, provided that the initial data is small enough.

By energy estimate, the Morawetz type estimate \eqref{8} and \eqref{2-decay2}, we see that \begin{eqnarray*}
  M_{k+1}&\le& C\ep+C \sum_{|\al|\le 2, |\be|\le 1}\|Z^\al\pa^\be F_p(u_k)\|_{L^1_t \dot H^{-1}}\\
  &\le &C \ep +C\sum_{|\al|\le 2}\|Z^\al F_p(u_k)\|_{L^1_t \dot H^{-1}}+C\sum_{|\al|\le 2 }\|Z^\al F_p(u_k)\|_{L^1_t L^2_x}\\
  &\le &C \ep +C\sum_{|\al|\le 2}\|\<x\>^{-1}Z^\al F_p(u_k)\|_{L^1_t L^1_{|x|}L^2_\theta}+C\sum_{|\al|\le 2 }\|Z^\al F_p(u_k)\|_{L^1_t L^2_x}\ .
\end{eqnarray*}
Then we can easily obtain the inequality $M_{k+1}\le C\ep +C M_k^p$, if we use the assumption \eqref{Fp-assume} on $F_p$, and the Sobolev embeddings $H^2_\theta\subset L^\infty_\theta$, $H^1_\theta\subset L^4_\theta$, $H^2_x\subset L^\infty_x$, $H^1_x\subset L^4_x$. Thus we can conclude that $M_k\le 2 C\ep$ for any $k\ge 0$.

\subsection{Adaption to the obstacle case: $n=4$ }
We shall consider wave equations on the exterior domain $\Omega\subset \R^n$ of a compact obstacle with metric $g$:
\begin{equation}\label{SLW}
\begin{cases}
\Box_\g u(t,x)=F_p(u(t,x)), \quad (t,x)\in \R_+\times \Omega
\\
u(0,\cdot)=f, \quad
\partial_t u(0,\cdot)=g,
\\
(Bu)(t,x)=0,\quad \text{on } \, \R_+\times \partial \Omega,
\end{cases}
\end{equation}
where for simplicity we take $B$ to either be the identity operator
(Dirichlet-wave equation) or the inward pointing normal derivative
$\partial_\nu$ (Neumann-wave equation). We shall assume that the nonlinear term behaves like $|u|^p$ when $u$ is small, and so we assume \eqref{Fp-assume} when $u$ is small, for $k=2$.

The operator $\Delta_\g$ is the Laplace-Beltrami operator associated with a
smooth, time independent Riemannian metric ${\mathbf g}_{jk}(x)$ which we
assume equals the Euclidean metric $\delta_{jk}$ for $|x|\ge R$, some $R$.
$\Omega=\R^n\backslash {\mathcal K}$ where ${\mathcal K}$ is a compact
subset of $|x|<R$ with smooth boundary.

In the work of \cite{DMSZ}, it shows that if $B=I$, then we have the existence result for $p>2$ as in the Strauss conjecture. For simplicity they have proved it with an extra assumption that the initial data $(f,g)$ are compactly supported. The idea of handling with a perturbed space is to consider the solution near and away from the obstacle separately. For that purpose we need KSS estimates, a variant of Hardy's inequality
\begin{equation}\label{h.2}
\|h\|_{\dot H^{-1}({\mathbb R}^n)} \le C\|h\|_{L^{2n/(n+2)}(|x|<1)} + C\|\, |x|^{-(n-2)/2} h\|_{L^1_{|x|} L^2_\theta(|x|>1)},
\end{equation}
and local energy decay (for $\Box_\g w=F$ with vanishing data)
\begin{equation}\label{2.13}
\sum_{|\alpha|\le 1}\| \partial^\alpha w\|_{L^2([0,T]\times \{x\in \Om: \, |x|<2\})} \le C\int_0^T\|F(t,\cdot)\|_{L^2(\Om)}\, dt.
\end{equation}
With these estimates we can have the Morawetz-KSS estimates with obstacle
\begin{align}\label{2.15}
&\sum_{|\alpha|\le N}\Bigl[
\bigl(\log(2+T)\bigr)^{-1/2}\|\langle x \rangle^{-1/2} Z^\alpha w\|_{L^2(S_T^\Om)}
+ \|\langle x\rangle^{-1/2-\delta}Z^\alpha w\|_{L^2(S_T^\Om)}
\\
&\qquad \qquad \qquad+\|Z^\alpha w(T,\cdot)\|_{L^2(\Om)} \Bigr],
\notag
\\
&\le C_\delta \int_0^T\sum_{|\alpha|\le N}\bigl(\, \|\,
\langle x\rangle^{-(n-2)/2}Z^\alpha  F(t,\cdot)\|_{L^1_{|x|} L^2_\theta(\Om)}
+\|Z^\alpha F(t,\cdot)\|_{L^2(\Om)}\, \bigr)\, dt\notag
\\
&\qquad \qquad \qquad+C\sum_{|\alpha|\le N-1}\Bigl[\, \sup_{0<s<T}\|Z^\alpha F(s,\cdot)\|_{L^2(\Om)}+\|Z^\alpha F\|_{L^2(S^\Om_T)}\Bigr]
\notag.
\end{align}
Where $S_T^\Om=[0,T]\times \Om$.
Now we can formally present the proof of existence results of \eqref{SLW}. We use a bump function $\eta\in C_0^\infty(\R)$, with $\eta(t)=1$ if $t<1/2$ and $\eta(t)=0$ if $t>1$. Set
$$u=\eta u+(1-\eta)u=u_0+w.$$
The local solution $u_0$ is assured by finite propagation of wave and compactly support property of data. For $w$, we set up the iteration as
\begin{multline*}M_k(T)=\sum_{|\alpha|\le 2}\Bigl(\|\langle x\rangle^{-1/2-\delta}Z^\alpha
w_k\|_{L^2(S_T^\Om)}+\sup_{0<t<T}\|Z^\alpha w_k(t,\cdot)\|_2\Bigr)
\\
+\sum_{|\alpha|\le 2}\Bigl(\|\langle x\rangle^{-1/2-\delta}Z^\alpha
w_k'\|_{L^2(S_T^\Om)}+\sup_{0<t<T}\|Z^\alpha w_k'(t,\cdot)\|_2\Bigr).
\end{multline*}
Now the boundedness of $M_k$ follows essentially from \eqref{2.15} and energy estimates. To technically make the iteration work, we just need the Sobolev inequalities
\eqref{2-decay}
and \begin{equation}\label{32}\|h\|_{L^4(\{|x|\in [R,R+1]\})}\le CR^{-3/4}\sum_{|\alpha|\le
1}\|Z^\alpha h\|_{L^2(\{|x|\in [R-1,R+2]\})}\ ,\end{equation}
where $h\in C^\infty(\R^4)$ and $R>1$.

\section{Weighted Strichartz estimates}

\subsection{A proof of the Strauss conjecture: $2\le n\le 4$}
As in the last section, we use continuity argument and iteration method to get the existence results for
        $$\Box u=|u|^p$$
with $p>p_c$ in dimensions $2$, $3$ and $4$.

The results can be found in  the works of Fang and Wang \cite{FaWa11} and Hidano, Metcalfe, Smith, Sogge and Zhou \cite{HMSSZ}. The main innovative idea is to utilize the following weighted Strichartz estimates (\eqref{3.5})
\begin{equation}
\Bigl\|\, |x|^{n/2-(n+1)/q-\gamma} e^{it|D|}\varphi \,
\Bigr\|_{L^q_tL^q_{|x|} L^2_\theta(\R_+\times \R^n)} \lesssim
\|\varphi\|_{\dot H^\gamma(\R^n)},
\end{equation}
for $2\le q\le \infty$ and $1/2-1/q<\gamma<n/2-1/q$. If we combine Duhamel formula and the duality of the trace estimate \eqref{7-g-Trace} with $s=1-\ga$, it is easy to obtain the following estimate which is the key ingredient of the proof of Strauss Conjecture:
\begin{multline}\label{Wstri}
\Bigl\|\, |x|^{n/2-(n+1)/p-\gamma} u \,
\Bigr\|_{L^p_tL^p_{|x|} L^2_\theta(\R_+\times \R^n)} \lesssim\\
\|f\|_{\dot H^\gamma(\R^n)}+\|g\|_{\dot H^{\gamma-1}(\R^n)}+\||x|^{-n/2+1-\gamma}F\|_{L^1_tL^1_{|x|}L^2_\om(\R_+\times\R^n)},
\end{multline}
where  $2\le p\le \infty$, $1/2-1/p<\gamma<n/2-1/p$, $1/2<1-\gamma<n/2$, and $u$ is the solution of the Cauchy problem
\begin{equation}\label{LWE}
\begin{cases}
\Box u(t,x)=F, \quad (t,x)\in \R_+\times \R^n
\\
u(0,\cdot)=f, \quad
\partial_t u(0,\cdot)=g.
\end{cases}
\end{equation}
Now if we set $\gamma=n/2-2/(p-1)$ so that $$-n/2+1-\gamma=p(n/2-(n+1)/p-\gamma)\ ,$$  the solution for $\Box u=|u|^p$ with data of size $\ep$ can be estimated as follows
\begin{eqnarray*}
  \||x|^{n/2-(n+1)/p-\gamma} u \|_{L^p_tL^p_{|x|}L^2_\theta}&\les& \|f\|_{\dot H^\ga}+\|g\|_{\dot H^{\ga-1}}+\||x|^{-n/2+1-\ga}|u|^p\|_{L^1_t \dot H^{-1}}\\
  &\les& \ep+\||x|^{-n/2+1-\ga}|u|^p\|_{L^1_t L^1_{|x|} L^2_\theta}\\
  &\les& \ep+\||x|^{n/2-(n+1)/p-\ga}|u|\|^p_{L^p_t L^p_{|x|} L^2_\theta}\les \ep+LHS^p\ .
\end{eqnarray*}

To give a rigorous proof, we need Sobolev embedding for $n=2,3,4$ to get around the last second inequality and set
$$M_k=\sum_{|\al|\leq 2}(\|Z^\al u\|_{L^\infty_t\dot H^\ga}+\||x|^{n/2-(n+1)/p-\gamma}Z^\al u \|_{L^p_tL^p_{|x|} L^2_\theta}).$$
Then the global existence result for $p>p_c$ are obtained by the classical energy estimates and standard continuity argument as in the last section. Here we notice that the requirement $p>p_c$ comes from the requirement $\ga>1/2-1/p$ and the observation \eqref{observ1}.

This method does not work for $n\geq 5$, due to the facts that $p_c(n)<2$ and our requirement for $p$ is $p\ge 2$.

\begin{rem}
If we use the endpoint trace estimate \eqref{7-Trace} instead of the trace estimates \eqref{7-g-Trace}, and the local in time KSS-type estimates \eqref{15}, we would get a local in time weighted-Strichartz estimate instead of \eqref{Wstri}. With that estimate, we can obtain the solution of \eqref{LW} with the sharp lifespan, for $p<p_c$ and $n=3$. See \cite{Yu11}.
\end{rem}

\subsection{Adaption to the obstacle case: $3\le n\le 4$}
When there exists an obstacle, new techniques are required due to reflection rays. By investigating the paramatrix of the wave near the obstacle (see for example \cite{MelSj}, \cite{Taylor76} and \cite{B}), people have proved several local energy decay estimates for the solution of the system. Making use of the local estimates, Smith and Sogge were able to get global Strichartz estimates for the wave equation with odd dimensions(\cite{SmSo00}, see also \cite{B} and \cite{M}). Stimulated by these ideas, in \cite{HMSSZ}, the authors created a mixed space $X$ and proved the Strauss conjecture for $n=3, 4$ with a nontrapping obstacle. The main idea is to use the Strichartz estimates in Minkowski space $\R\times\R^n$ for the solution away from the obstacle, and use local energy decay to control the part of the solution near the obstacle.

Precisely, the local energy decay is stated as follows,
\begin{equation}\label{localenergy}
\sum_{|\alpha|\le 1}\| \partial^\alpha u\|_{L^2(\R_+\times \{x\in \Om: \, |x|<R\})} \le C\left(\|f\|_{H^1(\Om)}+\|g\|_{L^2(\Om)}+\|F(t,\cdot)\|_{L^1_t L^2(\R_+\times\Om)}\right)\ .
\end{equation}
for the solution $u$ of the equation $\Box_\g u=F$ with data $(f,g)$,
where $f, g$ and $F(t,\cdot)$ are compactly supported in $|x|\le R$. With this estimate in hand,
by assuming the local-in-time abstract Strichartz estimates and the global-in-time Minkowski abstract Strichartz estimates, we get the following abstract Strichartz estimates for the system $\Box_\g u=0$ with data $(f,g)$
\beeq\label{global}\|u\|_{L^q_t X(\R\times\Om)}\les \|f\|_{\dot H^\gamma(\Om)}+\|g\|_{\dot H^{\ga-1}(\Om)}\ ,\eneq for $q>2$ and $\gamma\in [(3-n)/2, (n-1)/2]$.

To prove the Strauss conjecture under the setting of \eqref{SLW}, we define the space $X$ as follows
\beeq\label{X}\|u(t,\cdot)\|_X=\|u\|_{L^{s_\ga}(|x|<R)}+\||x|^{n/2-(n+1)/p-\ga}u\|_{L^p_{|x|} L^2_\theta(|x|>R)}\eneq
where $ s_\ga$ is such that $\dot H^\ga\subset L^{s_\ga}$. Then we apply Sobolev embedding to get the global Minkowski abstract Strichartz estimates
\begin{multline*}\|u\|_{L^p_t X(\R_+\times \R^n)}\les \|f\|_{\dot H^\ga}+\|g\|_{\dot H^{\ga-1}}+\\
\|F\|_{L^1_t L_x^{s'_{1-\ga}}(\R_+\times |x|<R)}+\||x|^{-n/2+1-\ga}F\|_{L^1_t L^1_r L^2_\theta(\R_+\times |x|>R)}.
\end{multline*}
The local-in-time abstract Strichartz estimates follow easily from the local energy decay estimate \eqref{localenergy} and the finite speed of propagation, and so is the global-in-time abstract Strichartz estimates \eqref{global}.

Now if we set
$$M_k = \sum_{|\alpha|\le 2}\,
\Bigl( \,
\bigl\|Z^\alpha u_k\bigr\|_{L^\infty_t\dot H^\gamma}+
\bigl\|\partial_tZ^\alpha u_k
\bigr\|_{L^\infty_t\dot H^{\gamma-1}}
+
\|Z^\al u\|_{L^p_tX}\Bigr ),
$$ the existence result for $n=3, 4$ follows from the same argument as in the Minkowski space.
The only reason for us to exclude the case $n=2$ is that we need to have $\ga\in [(3-n)/2, (n-1)/2]$ with $\ga=n/2-2/(p-1)$, for the abstract Strichartz estimates \eqref{global}.

\section{Generalized Strichartz estimates}
In this section, we present the generalized Strichartz estimates and illustrate the proof for the Strauss conjecture with $n=2$.
\subsection{Generalized Strichartz estimates}
The generalized Strichartz estimates for $n=2$ are stated as follows (we mainly follow \cite{SmSoWa10p}).
\begin{thm}\label{prop1.2}
Assume that
$(q,r)\ne (\infty,\infty)$
\begin{equation}\label{1.6}
q,r >2 \quad \text{and } \, \, \frac 1q < \frac 12 - \frac 1r\,,
\end{equation}
or $(q,r)=(\infty,2)$.
Then
\begin{equation}\label{1.7}
\bigl\|\, e^{it|D|}g \,  \bigr\|_{L^q_tL^r_{|x|}L^2_\theta({\mathbb R}\times {\mathbb R}^2)}
\le C_{q,r}\|g\|_{\dot H^{\gamma}({\mathbb R}^2)}\,,
\quad \gamma= 2(\tfrac 12-\tfrac 1r)-\tfrac 1q\,.
\end{equation}
\end{thm}
\begin{rem}
The radial estimates were first proven in \cite{FaWa06}. There are the endpoint and high dimensional versions of the generalized Strichartz estimates, see \cite{FaWa5} and \cite{JiWaYu10p}. For another type of estimate involving $L^q_t L^r_{x}$ norm and angular regularity, see \cite{St05}, \cite{FaWa11}.
\end{rem}

Since Hardy-Littlewood-Sobolev estimates give
$\dot H^{1-2/r}({\mathbb R}^2)\subset
L^r({\mathbb R}^2)$, $2\le r<\infty$, we clearly have
\begin{equation}\label{energy}
\bigl\|\, e^{i t |D|}f\, \bigr\|_{L^\infty_tL^r_{|x|}L^2_\theta}\le C_r\|f\|_{\dot H^{1-2/r}}\ .
\end{equation}
By scaling and Littlewood-Paley theory, it is easy to see that for the proof of \eqref{1.7}, we need only to prove
\begin{equation}\label{a}
\|e^{i t |D|}f\|_{L^q_tL^\infty_{|x|}L^2_\theta({\mathbb R}\times {\mathbb R}^2)}
\le C_q \|f\|_{L^2({\mathbb R}^2)}, \;\;
\text{if }q>2, \text{ and }
\Hat f(\xi)=0 \text{ if } |\xi|\notin [\tfrac 12,1].
\end{equation}

 By the support assumptions for
$\Hat f$ we have that
\begin{equation}\label{c}
\|f\|^2_{L^2({\mathbb R}^2)}\approx \int_0^\infty \int_0^{2\pi}|\Hat f(\rho(\cos\omega,\sin\omega))|^2
d\omega d\rho.
\end{equation}
We expand the angular part of $\Hat f$ using Fourier series and find that if
$\xi = \rho(\cos\omega,\sin\omega)$, then there are coefficients
$c_k(\rho)$ which vanish
for $\rho\notin [ 1/2,1]$, so that
$$\Hat f(\xi)=\sum_k c_k(\rho)\,e^{ik\omega}\,.$$
By \eqref{c} and Plancherel's theorem for $S^1$ and ${\mathbb R}$, we have
\begin{equation}\label{d}
\|f\|^2_{L^2({\mathbb R}^2)}
\approx \sum_k\int_{\mathbb R}|c_k(\rho)|^2\, d\rho \approx \sum_k\int_{{\mathbb R}}
|\Hat c_k(s)|^2 \, ds,
\end{equation}
where $\Hat c_k(s)$, $s\in {\mathbb R}$,
denotes the one-dimensional Fourier transform of $c_k(\rho)$.
Recall  that (see Stein and Weiss \cite{SW} p. 137)
\begin{equation}\label{e}
f(r(\cos\omega,\sin\omega))=(2\pi)^{-1}\sum_k
\Bigl(\, i^k \int_0^\infty J_{k}(r\rho)\,c_k(\rho)\,\rho \, d\rho\,
\Bigr) e^{ik\omega},
\end{equation}
where $J_k$, $k\in {\mathbb Z}$, is the $k$-th Bessel function, defined by
\begin{equation}\label{f}
J_k(y)=\frac{(-i)^k}{2\pi}\int_0^{2\pi} e^{iy\cos\theta -ik\theta}\, d\theta.
\end{equation}
By \eqref{e} and the support properties of $c_k$,
if we fix
$\beta\in C^\infty_0({\mathbb R})$ satisfying $\beta(\tau)=1$ for
$\tau\in[1/2,1]$ and $\beta(\tau)=0$ for $\tau\notin [1/4,2]$,
then with $\alpha(\rho)=\rho\,\beta(\rho)
\in {\mathcal S}({\mathbb R})$, we have
\begin{align*}
\bigl(e^{i t |D|}f\bigr)&(r(\cos\omega,\sin\omega))
\\
&=
(2\pi)^{-1}\sum_k\Bigl(\, i^k \int_0^\infty J_{k}(r\rho)\,e^{-it\rho}\,c_k(\rho)\,
\beta(\rho)\,\rho\, d\rho\, \Bigr)e^{ik\omega}
\\
&=(2\pi)^{-2}\sum_k\Bigl(\, i^k\int_0^\infty \int_{-\infty}^\infty J_{k}(r\rho)\,
e^{i\rho(s-t)}\,\Hat c_k(s)\,\alpha(\rho)\, ds\, d\rho\, \Bigr)e^{ik\omega}
\\
&=(2\pi)^{-3}\sum_k\Bigl(\,
\int_0^\infty\int_{-\infty}^\infty \int_0^{2\pi}
e^{i\rho r\cos\theta}e^{-ik\theta}e^{i\rho(s-t)}\,\Hat c_k(s)\,
\alpha(\rho)\, d\theta\, ds \,d\rho\, \Bigr)e^{ik\omega}
\\
&=(2\pi)^{-3}\sum_k\Bigl(\,  \int_{-\infty}^\infty \int_0^{2\pi}
e^{-ik\theta}\Hat \alpha\bigl((t-s)-r\cos\theta \bigr)\, \Hat c_k(s) \,
d\theta \, ds \, \Bigr)e^{ik\omega}.
\end{align*}
As a result, we have that for any $r\ge0$,
\begin{multline}\label{g}
\int_0^{2\pi}\Bigl| \, \bigl(e^{i t |D|}f\bigr)
(r(\cos\omega,\sin\omega))\, \Bigr|^2 \, d\omega
\\
=(2\pi)^{-5} \sum_k\,\Bigl|\, \int_{-\infty}^\infty \int_0^{2\pi} e^{-ik\theta}\,
\Hat \alpha\bigl((t-s)-r\cos\theta\bigr)\, \Hat c_k(s) \,d\theta\,ds\, \Bigr|^2.
\end{multline}

To estimate the right side we shall use the following Lemma.
\begin{lem}[Lemma 2.1 of \cite{SmSoWa10p}]\label{mainest}
  Let $\alpha\in {\mathcal S}({\mathbb R})$ and $N\in {\mathbb N}$ be fixed.  Then there is a uniform constant $C$, which is independent of $m\in {\mathbb R}$ and $r\ge0$, so that the following inequalities hold.  First,
$$
\int_0^{2\pi} |\alpha(m-r\cos\theta)|\, d\theta \le C \langle \, m\, \rangle^{-N}, \quad
\text{if } \, \, 0\le r\le 1, \, \, \text{or }\, \, |m|\ge 2r.
$$
If $r>1$ and $|m|\le 2r$ then
$$\int_0^{2\pi}|\alpha(m-r\cos\theta)|\, d\theta
\le C\Bigl(\, r^{-1}+r^{-\frac12}\langle \, r-|m|\, \rangle^{-\frac12}\, \Bigr).
$$Consequently, if $\delta>0$, there is a constant $A_\delta$, which is independent of $t\in {\mathbb R}$ and $r\ge 0$ so that
\begin{equation}\label{ssw3}
\int_{-\infty}^\infty\left( \, \int_0^{2\pi} \langle \, t-s \, \rangle^{\frac12-\delta}\,
|\alpha((t-s)-r\cos\theta)|\, d\theta \, \right)^2 \, ds \le A_\delta.
\end{equation}
\end{lem}

If we apply \eqref{ssw3} and \eqref{g} along with the Schwarz inequality, we
conclude that if $f$ is as in \eqref{a}, then for $\delta>0$
$$\Bigl\| \, e^{i t |D|}f\, \Bigr\|_{L^\infty_{|x|}L^2_\theta}^2
\le B_\delta \sum_k\int_{-\infty}^\infty \, \bigl|\, \langle\, t-s\,
\rangle^{-\hf+\delta}
\Hat c_k(s)\, \bigr|^2 \, ds,$$
which, by Minkowski's inequality and \eqref{d}, in turn yields \eqref{a}.

\subsection{A new proof of the Strauss conjecture for $n=2$}\label{sec-Glassey}
In this subsection, we illustrate how Theorem~\ref{prop1.2}
implies estimates that can be used to prove
Glassey's existence theorem \cite{G} for $\square u =|u|^p$ when $n=2$.
Specifically, if $u$
solves the wave equation \eqref{LWE} for ${\mathbb R}\times {\mathbb R}^2$,
then
\begin{equation}
\label{ssw15}
\|u\|_{L^q_tL^r_{|x|}L^2_\theta}+\|u\|_{L^\infty_t\dot H^\gamma}\lesssim
\|f\|_{\dot H^\gamma}+\|g\|_{\dot H^{\gamma-1}}+\|F\|_{L_t^{\tilde q'}L^{\tilde r'}_{|x|}L^2_\theta},
\end{equation}
assuming that $q,r, \tilde q, \tilde r>2$ with
$(q,r), (\tilde q,\tilde r) \ne (\infty, \infty)\,,$
$ 1/q<1/2- 1/r$, $1/{\tilde q}< 1/2 -  1/{\tilde r}$, and
\begin{equation}\label{ssw16}
\gamma=1-\frac 2r-\frac 1q\,, \quad \text{and } \,
1-\gamma = 1-\frac 2{\tilde r}-\frac 1{\tilde q}\,.
\end{equation}
 Clearly, \eqref{ssw15} follows from \eqref{1.7} and
energy estimates if the forcing term, $F$, in \eqref{LWE} vanishes.  Since we
are
assuming \eqref{ssw16} and since $\tilde q'<q$, the estimates for the
inhomogeneous wave
equation follow from an application of the Christ-Kiselev lemma \cite{ChKi01}
(cf. \cite{So08},
pp. 136--141).

Here, let us just give the key estimate for the new proof of the Strauss conjecture with $n=2$.
Considering the subconformal range $(3+\sqrt{17})/2=p_c<p<5$, as a special case of \eqref{ssw15},
we have \beeq\label{ssw20}
\| u\|_{L^{\frac{(p-1)p}{2}}_tL^p_{|x|}L^2_\theta}+\| u\|_{L^\infty_t\dot
H^{\gamma_p}}
\lesssim
\| f\|_{\dot H^{\gamma_p}}+\| g\|_{\dot H^{\gamma_p-1}}
+\| F\|_{L^{\frac{p-1}{2}}_t L^1_{|x|}L^2_\theta}\,.
\eneq
The temporary assumption that $p<5$ is needed to ensure that $(p-1)/2<2$,
and, therefore,
$[(p-1)/2]'>2$, which is the first part of the assumptions for \eqref{ssw15}.
The more serious
assumption that $p>p_c$  is equivalent to the second part of
\eqref{1.6} for the exponents on the left side of \eqref{ssw20}. That is,
for $p>0$,
$$\frac2{p(p-1)}<\frac12-\frac1p \, \iff \, p>p_c.$$

Using \eqref{ssw20}, we can easily  solve \eqref{LW} by an
iteration argument
for $p_c<p<5$, provided that the initial data is small.

\subsection{Adaption to the obstacle case}

It is natural to try to extend the previous proof in the obstacle setting, following the argument of \cite{HMSSZ}. Precisely, mimicking the $X$ space as in \eqref{X},
we define $X$ space as follows
$$\|h\|_{X_{r,\gamma}}= \|h\|_{L^{s_\gamma}(|x|<3R)}+\|h\|_{L^r_{|x|}L^2_\theta(|x|>2R)},
\quad \text{with } \, \gamma=1-\frac 2{s_\gamma}\,.$$
By local energy estimates (\cite{V}, \cite{B}) and the abstract Strichartz estimates of \cite{HMSSZ}, we know that we have global estimate \eqref{global}, if $q>2$ and $\gamma=1/2$. This estimate is obviously not enough for the Strauss conjecture.

As noted in \cite{SmSoWa10p}, even though we can only directly prove Strichartz estimates
involving Sobolev regularity of $1/2$, for some applications if we interpolate
with trivial (energy) estimates, this is enough.
In this case, by interpolation, we can prove the global estimate \eqref{global} for $q>2$ and $0<\gamma<1$. Then, despite some technical difficulties with obstacles, we can in principle adapt the proof in Section \ref{sec-Glassey} to the problem with obstacles.

\section{Appendix}
 \subsection{Hardy's inequality}
Let us first record a fundamental inequality.
\begin{thm}[Hardy's inequality]
Let $0\le s<n/2$,  we have
\beeq\left\|\frac{u}{|x|^s}\right\|_{L^2}\le C \|u\|_{\dot H^{s}}.\eneq
\end{thm}
See Tao \cite{Tao06} (Lemma A.2 p334) for a simple proof by using Littlewood-Paley decomposition.

Here, we present an elementary proof for $s=1$ by using integration by parts: when $s=1$ and $n\ge 3$, we have
\beeq \left\|\frac{u}{|x|}\right\|_{L^2} \le
\frac{2}{n-2} \| \nabla u\|_{L^2}. \eneq

In fact, we have \begin{eqnarray*} \left\|\frac{u}{|x|}\right\|^2_{L^2} &=& \int_{S^{n-1}}
\int_0^\infty |u|^2
r^{n-3} dr d\omega\\
& =& \int_{S^{n-1}} \int_0^\infty |u|^2 \pa_r \frac{r^{n-2}}{n-2}
dr d\omega \\
&\le &  \frac{2}{n-2} \int_{S^{n-1}} \int_0^\infty |\pa_r u|
|u| r^{n-2} d r d \omega \\
&\le &  \frac{2}{n-2} \| |\pa_r u|
|u|  r^{-1} \|_{L^1} \\
&\le &  \frac{2}{n-2} \||x|^{-1} u\|_{L^2}\| \pa_r
u  \|_{L^2} .\end{eqnarray*}
As far as we know, the proof of this type has occurred in the works of Klainerman, Li, Sideris, Yu, Zhou and many others, since, as early as 1980's.

\subsection{Trace lemma}

  One version of the trace lemma can be stated as follows
$$r^{\frac{n-1}{2}} \|f(r \omega)\|_{L^2_{\omega}}
  \les \|   f\|_{\dot{B}^{\frac{1}{2}}_{2,1}}\ .$$
Moreover, for $s\in (1/2,n/2)$, we have $$r^{\frac{n}{2}-s}
\|f(r\omega)\|_{L^2_\omega}
  \les \|f\|_{\dot H^s_x}\ .$$ For more delicate version and proof, see Fang and Wang \cite{FaWa11}.

Here, we want to give an elementary proof of these two estimates.
Our proof is inspired by the paper of Li and Yu \cite{LiYu91}

\begin{lem}\label{lem-trace}
  For $s\in [1, n/2)$, we have
$$r^{\frac{n}{2}-s}
\|f(r\omega)\|_{L^2_\omega}
  \les \|f\|_{\dot H^s_x}\ .$$
Moreover, for any $n\ge 1$,
\beeq\label{eq-trace}r^{\frac{n-1}{2}} \|f(r
\omega)\|_{L^2_{\omega}}
  \les \|   f\|_{L^2}^{1/2}\|   f\|_{\dot H^1}^{1/2}\ .\eneq
\end{lem}

For the second inequality,
\begin{eqnarray*}
\int_{\Sp^{n-1}}|f (r\omega)|^2
d\omega&=&-\int_{\Sp^{n-1}}\int_r^\infty \pa_\la|f (\la\omega)|^2  d\la d\omega\\
&\le & 2\ r^{1-n} \int_{\Sp^{n-1}}\int_r^\infty |f| |\pa_\la f|
\la^{n-1}d\la d\omega\\
&\le & 2\ r^{1-n} \|f\|_{L^2} \|\pa_\la f\|_{L^2}
\end{eqnarray*}

For the first inequality, since $s\ge 1$,
\begin{eqnarray*}
\int_{\Sp^{n-1}}|f (r\omega)|^2
d\omega&=&-\int_{\Sp^{n-1}}\int_r^\infty \pa_\la|f (\la\omega)|^2  d\la d\omega\\
&\le & 2\ r^{2s-n} \int_{\Sp^{n-1}}\int_r^\infty |f| |\pa_\la f|
\la^{n-2s}d\la d\omega\\
&\le & 2\ r^{2s-n} \int_{\Sp^{n-1}}\int_r^\infty |\la^{-s}f|
|\la^{1-s} \pa_\la f|
\la^{n-1}d\la d\omega\\
&\le & 2\ r^{2s-n} \|\la^{-s} f\|_{L^2} \|\la^{1-s} \pa_\la
f\|_{L^2}\\
&\le& C r^{2s-n} \|f\|_{\dot H^s}^2\ ,
\end{eqnarray*} where in the last inequality we have used the
Hardy's inequality.

As we can see, Lemma \ref{lem-trace} is enough for the proof of the general trace estimates.
\begin{thm}[Trace lemma]
  For $s\in (1/2, n/2)$, we have
\beeq\label{7-g-Trace}r^{\frac{n}{2}-s} \|f(r\omega)\|_{L^2_\omega}
  \les \|f\|_{\dot H^s_x}\eneq
Moreover, for any $n\ge 1$, \beeq\label{7-Trace}r^{\frac{n-1}{2}}
\|f(r \omega)\|_{L^2_{\omega}}
  \les \|   f\|_{\dot B^{1/2}_{2,1}}\eneq
\end{thm}
\begin{prf}
 Applying \eqref{eq-trace} to the Littlewood-Paley projection $P_\la f$ with frequency of size $\la$, we see that
$$ r^{\frac{n-1}{2}} \|P_\la f(r
\omega)\|_{L^2_{\omega}}
  \les \la^{1/2}\|   P_\la f\|_{L^2}\ .$$
By Littlewood-Paley decomposition, we arrived
$$ r^{\frac{n-1}{2}} \|f(r
\omega)\|_{L^2_{\omega}}
  \les \|f\|_{\dot B^{1/2}_{2,1}}\ .$$

We are remained to prove the first inequality for $1/2<s<1$. Once again, for fixed $s$, with $\theta=2s-1$,
we have
$$r^{\frac{n}{2}-s}
\|f(r\omega)\|_{L^2_\omega}
  \les (r^{\frac{n-1}{2}}
\|f(r\omega)\|_{L^2_\omega}
 )^{1-\theta}(r^{\frac{n-2}{2}}
\|f(r\omega)\|_{L^2_\omega}
 )^{\theta}\les  \|f\|_{\dot B^{1/2}_{2,1}}^{1-\theta}\|f\|_{\dot H^1}^\theta \ .$$
Applying this estimate to $P_\la f$, we get
$$r^{\frac{n}{2}-s}
\|P_\la f(r\omega)\|_{L^2_\omega}
  \les \la^s \|P_\la f\|_{L^2}\ ,$$
and so
$$ r^{\frac{n}{2}-s} \|f(r
\omega)\|_{L^2_{\omega}}
  \les \|f\|_{\dot B^{s}_{2,1}}\ .$$

Recall that for $s_0\neq s_1$, we have the following fact of the real interpolation (see \cite{BeLo} p152 (1))
$$[\dot B^{s_1}_{2,1}, \dot B^{s_2}_{2,1}]_{\theta,2}=\dot H^{s}\ ,$$
which tells us that we can actually have
$$ r^{\frac{n}{2}-s} \|f(r
\omega)\|_{L^2_{\omega}}
  \les \|f\|_{\dot H^{s}}\ $$
for $1/2<s<1$.
\end{prf}

\subsection{Morawetz type estimates}
With the help of the trace estimates, we can prove the Morawetz type estimates.
\begin{thm}[Morawetz type estimates]
For $s\in (1/2,n/2)$, we have
\begin{equation}\label{3.4}
\bigl\|\, |x|^{-s} e^{it|D|}\varphi \|_{L^2(\R_+ \times \R^n)}
\lesssim \|\varphi\|_{\dot H^{s-\frac12}\R^n)}, \quad
\frac12<s<\frac{n}2\,.
\end{equation}
Moreover, for any $n\ge 1$, we have
\begin{equation}\label{3.6}
\bigl\|\, e^{it|D|}\varphi \|_{L^2(\R_+ \times \{|x|\sim 1\})}
\lesssim \|\, \varphi\|_{L^2(\R^n)} \,.
\end{equation}
\end{thm}
As a corollary, we have $$\sup_{R} R^{-1/2}\bigl\|\, e^{it|D|}\varphi
\|_{L^2(\R_+ \times \{|x|\les R\})} \lesssim \|\,
\varphi\|_{L^2(\R^n)} \,.$$

The inequality \eqref{3.6} is also known as local energy estimate, which dates back to the works of Morawetz, Strauss and others in 1970's. A more general version has occurred in Lemma 2.2 of \cite{SmSo00}. The inequality \eqref{3.4} is sometimes called generalized Morawetz estimate, see \cite{Ho97}, \cite{HMSSZ}, and \cite{FaWa11}.

\begin{rem}\label{rq91}
  As a consequence of \eqref{3.6}, we have
  \beeq\label{8}
\| \<x\>^{-b} e^{i t |D|} f\|_{L^2_{t, x} } \les \| f\|_{L^2_x}
\eneq for any $b>1/2$, and
  \beeq\label{9}
\| |x|^{-b} e^{i t |D|} f\|_{L^2_{t, |x| \le 1} } \les \| f\|_{L^2_x}
\eneq for any $b<1/2$.
\end{rem}

Note that by applying \eqref{7-g-Trace} to the Fourier transform of
$v$, we see that it is equivalent to the uniform bounds
$$
\Bigl(\, \int_{S^{n-1}}|\Hat v(\lambda\omega)|^2\, d\omega\,
\Bigr)^{1/2}\lesssim \lambda^{-\frac{n}2+s}\|\, |x|^s
v\|_{L^2(\R^n)}\,, \quad \lambda>0\,, \quad \frac12<s<\frac{n}2\,,
$$
which by duality is equivalent to
\begin{equation}\label{3.3}
\Bigl\|\, |x|^{-s}\int_{S^{n-1}}h(\omega) e^{i\lambda x\cdot
\omega}\, d \omega \,  \Bigr\|_{L^2_x(\R^n)}\lesssim
\lambda^{s-\frac{n}2}\|h\|_{L^2_\omega(S^{n-1})}\,,
\end{equation}
for $\lambda>0$ and fixed $1/2<s<n/2$.  Using this estimate we can
obtain \eqref{3.4}.

In fact, recall that $$ \mathfrak{F}_{tx} e^{it|D|}\varphi=\int_{\R}
e^{-it(\tau-|\xi|)}\hat \varphi(\xi)dt \sim \delta(\tau-|\xi|)\hat
\varphi(\xi)\ ,$$
$$ \mathfrak{F}_{t} e^{it|D|}\varphi\sim
\int_{\R^n}e^{ix\cdot \xi}\delta(\tau-|\xi|)\hat\varphi(\xi)d\xi
=\int_{S^{n-1}}e^{i\tau x\cdot \omega}\tau^{n-1}\hat\varphi(\tau
\omega)d\omega\ .$$ Thus, by Plancherel's theorem with respect
to the $t$-variable, we find that the square of the left side of
\eqref{3.4} equals
\begin{multline*}
(2\pi)^{-1} \int_0^\infty \int_{\R^n}\Bigl| \,
|x|^{-s}\int_{S^{n-1}} e^{ix\cdot \rho \omega}\rho^{n-1}\Hat
\varphi(\rho\omega)\, d \omega\, \Bigr|^2\, dx \, d\rho
\\
\lesssim \int_0^\infty \int_{S^{n-1}}\rho^{2(n-1)}| \Hat
\varphi(\rho\omega)|^2 \, \rho^{2s-n}\, d\omega d\rho = \|\,
|D|^{s-\frac12}\varphi\|_{L^2(\R^n)}^2,
\end{multline*}
by using \eqref{3.3} in the first step.

If we apply \eqref{7-Trace} instead, we see that for $s=1/2$
$$
\Bigl(\, \int_{S^{n-1}}|\Hat v(\lambda\omega)|^2\, d\omega\,
\Bigr)^{1/2}\lesssim \lambda^{\frac{1-n}2}\|\, |x|^{1/2} \phi(x
2^{-j}) v\|_{l_j^1 L^2(\R^n)}\,, \quad \lambda>0\,,
$$
which by duality is equivalent to
\begin{equation}\label{3.5}
\Bigl\|\, |x|^{-1/2}\phi(x 2^{-j})\int_{S^{n-1}}h(\omega)
e^{i\lambda x\cdot \omega}\, d \omega \, \Bigr\|_{l_j^\infty
L^2_x(\R^n)}\lesssim
\lambda^{\frac{1-n}2}\|h\|_{L^2_\omega(S^{n-1})}\,,
\end{equation}
for $\lambda>0$ and fixed $1/2<s<n/2$.  Using this estimate we can
obtain \eqref{3.6}.

In fact, recall that
$$ \mathfrak{F}_{t} e^{it|D|}\varphi\sim
\int_{\R^n}e^{ix\cdot \xi}\delta(\tau-|\xi|)\hat\varphi(\xi)d\xi
=\int_{S^{n-1}}e^{i\tau x\cdot \omega}\tau^{n-1}\hat\varphi(\tau
\omega)d\omega\ .$$ Thus, by after Plancherel's theorem with respect
to the $t$-variable, we find that the square of the left side of
\eqref{3.6} equals
\begin{multline*}
(2\pi)^{-1} \int_0^\infty \int_{|x|\sim 1}\Bigl| \,
|x|^{-s}\int_{S^{n-1}} e^{ix\cdot \rho \omega}\rho^{n-1}\Hat
\varphi(\rho\omega)\, d \omega \, \Bigr|^2\, dx \, d\rho
\\
\lesssim \int_0^\infty \int_{S^{n-1}}\rho^{2(n-1)}| \Hat
\varphi(\rho\omega)|^2 \, \rho^{1-n}\, d \omega d\rho = \|\,
\varphi\|_{L^2(\R^n)}^2,
\end{multline*}
by using \eqref{3.5} in the first step.

\subsection[Morawetz-KSS estimates]{Morawetz-KSS estimates}

Given the Morawetz type estimates, we can combine with the
energy estimate to obtain the Keel-Smith-Sogge type estimates (KSS in short).
\begin{thm}[Morawetz-KSS estimates]
Let $n\ge 1$. For any $T>0$, we have \beeq\label{eq-KSS}\log(2+T)^{-\frac{1}{2}}\|
\langle x\rangle^{-\frac{1}{2}} e^{i t |D|} f\|_{L_{[0,T]}^2
L^2_{x}}\les \|f\|_{L^2_x}.\eneq Moreover, for any $\ep\in
(0,\frac{1}{2}]$ and $\delta>0$, then \beeq\label{13}\| \langle
x\rangle^{-\ep}|x|^{\ep-\frac{1}{2}} e^{i t |D|} f\|_{L_{[0,T]}^2
L^2_{x}}\les \min(\log(2+T)^{\frac{1}{2}}, T^\ep) \| f\|_{L^2_x}.
\eneq \beeq\label{14}\| \langle
x\rangle^{-\ep-\delta}|x|^{\ep-\frac{1}{2}} e^{i t |D|} f\|_{L_{\R}^2
L^2_{x}}\les \| f\|_{L^2_x}.\eneq
\end{thm}
The estimate \eqref{eq-KSS} was first proven in Keel, Smith and Sogge \cite{KeSmSo02} for $n=3$. The estimates of this type were developed drastically afterwards (see e.g. \cite{Met04-1}, \cite{HiYo05}, \cite{MetSo06}, \cite{Al06}, \cite{SW} and \cite{HiWaYo}). The general versions \eqref{13} and \eqref{14} were proven in Hidano and Yokoyama \cite{HiYo05}, which were also generalized to the situation of the wave equation with variable coefficients in \cite{HiWaYo}.

\begin{prf}
At first, the estimate \eqref{14} is clear from Remark \ref{rq91}. We need only to give the proof of \eqref{13}.

By Plancherel's theorem, we have
$$\|e^{it |D|}f\|_{L^2_x}=\|f\|_{L^2_x},$$
thus $$T^{-1/2}\|e^{it |D|}f\|_{L^2_{t\in [0,T]}
L^2_x}=\|f\|_{L^2_x}.$$ Recall also the Morawetz type estimate
\eqref{3.6},
$$ R^{-\frac{1}{2}} \| e^{i t |D|} f\|_{L^2_{t, x\in B(0,
R) } } \les \| f\|_{L^2_x}.$$ We consider first the case when $T\ge
1$. If $|x|\le  1$, then since $\ep>0$, we have
$$\| \langle
x\rangle^{-\ep}|x|^{\ep-\frac{1}{2}} e^{i t |D|} f\|_{L_{[0,T]}^2
L^2_{|x|\le 1}} \les \sum_{j\le 0}\||x|^{\ep-\frac{1}{2}} e^{i t |D|}
f\|_{L_{T}^2 L^2_{|x|\simeq 2^j}} \les \sum_{j\le 0}2^{j\ep} \|
 f\|_{L^2_x}\les \|
f\|_{L^2_x}.$$ For $1\le |x|\le T$, we have
$$\| \langle
x\rangle^{-\ep}|x|^{\ep-\frac{1}{2}} e^{i t |D|} f\|_{L_{[0,T]}^2
L^2_{|x|\in [1,T]}}^2 \les \sum_{ 0\le j \le \ln(2+T)
}\||x|^{-\frac{1}{2}} e^{i t |D|} f\|^2_{L_{T}^2 L^2_{|x|\simeq
2^j}}\les  \ln(2+T)\|
 f\|_{L^2_x}^2.$$ For the remained $ |x|\ge T$, we have
$$\| \langle
x\rangle^{-\ep}|x|^{\ep-\frac{1}{2}} e^{i t |D|} f\|_{L_{[0,T]}^2
L^2_{|x|\ge T}} \les T^{-1/2}\| e^{i t |D|} f\|_{L_{T}^2 L^2_{x}}
\les \|
 f\|_{L^2_x}.$$ This complete the proof for the case $T\ge 1$. For
 $T\le 1$, we consider three cases. If $|x|\le  T$, then
$$\| \langle
x\rangle^{-\ep}|x|^{\ep-\frac{1}{2}} e^{i t |D|} f\|_{L_{[0,T]}^2
L^2_{|x|\le T}} \les \sum_{j\le \ln T}\||x|^{\ep-\frac{1}{2}} e^{i t
|D|} f\|_{L_{T}^2 L^2_{|x|\simeq 2^j}} \les \sum_{j\le \ln T}
2^{j\ep} \|
 f\|_{L^2_x}\les T^\ep \|
f\|_{L^2_x}.$$ For $T\le |x|\le 1$, we have
$$\| \langle
x\rangle^{-\ep}|x|^{\ep-\frac{1}{2}} e^{i t |D|} f\|_{L_{[0,T]}^2
L^2_{|x|\in [T,1]}} \les T^{\ep-1/2}\|e^{i t |D|} f\|_{L_{T}^2
L^2_{x}}\les  T^\ep \|
 f\|_{L^2_x}.$$ For $ |x|\ge 1$ and if $\ep\le 1/2$, we have
$$\| \langle
x\rangle^{-\ep}|x|^{\ep-\frac{1}{2}} e^{i t |D|} f\|_{L_{[0,T]}^2
L^2_{|x|\ge 1}} \les \|e^{i t |D|} f\|_{L_{T}^2 L^2_{x}}\les
T^{\frac{1}{2}} \|  f\|_{L^2_x} \les T^{\ep} \|
 f\|_{L^2_x}.$$
This completes the proof of the KSS estimate.\end{prf}

\begin{rem}
By similar arguments, we can get the following
\beeq\label{15}T^{-\ep}\| |x|^{\ep-\hf} e^{i t |D|}
f\|_{L_{[0,T]}^2 L^2_{x}}\les \|f\|_{L^2_x}\eneq for $\ep\in (0,
\frac{1}{2}]$. See Hidano \cite{Hi07}.\end{rem}

In fact, the estimate is scale invariant. Without loss of generality, we can let $T=1$. For the part of
$|x|\le 1$, it is obvious from Remark \ref{rq91}. Else, if
$|x|\ge 1$, the estimate is weaker than the energy estimate,
$$\| |x|^{\ep-\hf} e^{i t |D|} f\|_{L_{[0,1]}^2 L^2_{|x|\ge 1}}\les \| e^{i t |D|} f\|_{L_{[0,1]}^2 L^2_{x}}
\les \|f\|_{L^2_x}\ .$$

\subsection{Weighted Strichartz estimates}
If we interpolate between \eqref{7-g-Trace} and \eqref{3.4} we
conclude that, for $2\le q\le \infty$,
\begin{equation}\label{3.5}
\Bigl\|\, |x|^{\frac{n}2-\frac{n+1}q-\gamma} e^{it|D|}\varphi \,
\Bigr\|_{L^q_tL^q_{|x|} L^2_\theta(\R_+\times \R^n)} \lesssim
\|\varphi\|_{\dot H^\gamma(\R^n)}, \quad
\frac12-\frac1q<\gamma<\frac{n}2-\frac1q.
\end{equation}
The estimates of this type were first proved by Hidano \cite{Hi07} in the radial case, and then generalized in Fang and Wang \cite{FaWa11} and Hidano, Metcalfe, Smith, Sogge and Zhou \cite{HMSSZ}. For more general estimates involving mixed norm in $t$ and $|x|$, see our recent work \cite{JiWaYu10p}.

\end{document}